\documentclass[final,nomarks]{dmtcs-episciences}

\usepackage[utf8]{inputenc}
\usepackage{subfigure}

\usepackage{amssymb}
\usepackage{amsmath}
\usepackage{amsthm,MnSymbol}

\newtheoremstyle{meiner} 
    {10pt}{5pt}          
    {\itshape}           
    {}                   
    {\sffamily\bfseries} 
    {}                   
    { }                  
    {}                   
\theoremstyle{meiner}

\newtheorem{theorem}{Theorem}
\newtheorem{lemma}{Lemma}
\newtheorem{proposition}{Proposition}

\def\Claim#1.{\par\smallskip{\ni\bfseries Claim #1.}}
\def\Fact#1.{\par\smallskip{\ni\bfseries Fact #1.}}

\usepackage{bibitemsep}

\def\ni{\noindent}

\def\Proof{\par{\ni\bfseries Proof. }}
\def\qed{\hfill\fbox{\hbox{}}\medskip}
\def\qedclaim{\hfill$\triangle$\smallskip}

\def\term#1{{\em #1}\marginpar{{\setlength\baselineskip{6pt}\raggedright\tiny\em #1\par}}}
\def\nct#1{{\em #1}}

\def\ITEMMACRO #1 ??? #2 ???{\par\medskip\noindent%
\hangindent=#2em\setbox0\hbox{#1 \kern5pt}%
\ifdim\wd0<\hangindent\setbox0\hbox to\hangindent{\hss#1\quad}\fi%
\box0\ignorespaces}

\def\Item#1{\ITEMMACRO #1 ??? 2.5 ???}

\def\Bitem{\Item{\hss$\bullet$}}

\input psfig.sty
\usepackage{wrapfig}

\DeclareMathOperator{\E}{\mathbb{E}}
\def\NN{\mathbb{N}}
\DeclareMathOperator{\outdeg}{{\mathrm outdeg}}

\def\hf{\hat{f}}
\def\cL{{\cal L}}
\def\Tmix{\tau_\text{mix}}
\def\vX{\vec{X}}
\def\vY{\vec{Y}}
\def\vG{\vec{G}}
\def\Goct{G_\textsf{oct}}
\def\Gmin{\vec{G}_\textsf{min}}

\author{Stefan Felsner\thanks{Partially supported by DFG FE-340/11-1}
  \and Daniel Heldt}
\title[Markov Chains on Degree Constrained Orientations]%
 {Mixing Times of Markov Chains on Degree Constrained
  Orientations of Planar Graphs}
\affiliation{
  Institut f\"ur Mathematik,
  Technische Universit\"at Berlin}

\received{2016-2-10}
\accepted{2016-12-5}
\publicationdetails{18}{2016}{3}{20}{1376}
\begin{document}
\publicationdetails{18}{2016}{3}{20}{1376}
\maketitle
\begin{abstract}
We study Markov chains for $\alpha$-orientations of plane graphs,
these are orientations where the outdegree of each vertex is 
prescribed by the value of a given function $\alpha$. The set of
$\alpha$-orientations of a plane graph has a natural distributive
lattice structure. The moves of the up-down Markov chain on this distributive
lattice corresponds to reversals of directed facial cycles in the 
$\alpha$-orientation. We have a positive and several negative results regarding
the mixing time of such Markov chains.

A 2-orientation of a plane quadrangulation is an orientation where
every inner vertex has outdegree 2. We show that there is a class of plane quadrangulations
such that the up-down Markov chain on
the 2-orientations of these quadrangulations is slowly mixing. 
On the other hand the chain is rapidly mixing on 2-orientations of 
quadrangulations with maximum degree at most 4.

Regarding examples for slow mixing we also revisit the case of
3-orientations of triangulations which has been studied before by
Miracle et al. Our examples for slow mixing are
simpler and have a smaller maximum degree, Finally we present the
first example of a function~$\alpha$ and a class of plane
triangulations of constant maximum degree such that the up-down Markov
chain on the $\alpha$-orientations of these graphs is slowly mixing.
\end{abstract}

\section{Introduction}

Let $G=(V,E)$ be a graph and let $\alpha:V\to \NN$ be a function, an
$\alpha$-orientation of $G$ is an orientation with $\outdeg(v) = \alpha(v)$
for all vertices $v\in V$.  A variety of interesting combinatorial structures
on planar graphs can be modeled as $\alpha$-orientations. Examples are
spanning trees, Eulerian orientations, Schnyder woods of triangulations,
separating decompositions of quadrangulations. These and further examples are
discussed in~\cite{Fe-04} and \cite{FZ-08}. In this paper we are interested in
Markov chains to sample uniformly from the $\alpha$-orientations of a given
planar graph $G$ for a fixed $\alpha$. 

A uniform sampler may be used to get data for a statistical approach
to typical properties of $\alpha$-orientations. Under certain
conditions such a chain can be used for approximate counting of
$\alpha$-orientations.  Counting $\alpha$-orientations is \#P-complete
in general.  Mihail and Winkler~\cite{MW-96} have shown that counting
Eulerian orientations is \#P-complete.  Creed~\cite{Cr-09} has shown
that this counting problem remains \#P-complete on planar
graphs. Further examples of \#P-complete variants of counting
$\alpha$-orientations are given in~\cite{FZ-08}.  In~\cite[Section
6.2]{FZ-08} it is shown that counting $\alpha$-orientations can be
reduced to counting perfect matchings of a related bipartite
graph. The latter problem can be approximately solved using the
celebrated algorithm of Jerrum, Sinclair and
Vigoda~\cite{JSV-Permanent-04} or its improved version by Bez\'akov\'a
et al.~\cite{BSVV-08}. These algorithms build on random sampling.

For sampling $\alpha$-orientations of plane graphs, however, there is
a more natural Markov chain. The reversal of the orientation of a
directed cycle in an $\alpha$-orientation yields another
$\alpha$-orientation. If $G$ is a plane graph and $\vec{G},\vec{G}'$
are $\alpha$-orientations of $G$, then we define $\vec{G} < \vec{G}'$
whenever $\vec{G}'$ is obtained by reverting a clockwise cycle
of $\vec{G}$.  In~\cite{Fe-04} it has been shown that this order
relation makes the set of $\alpha$-orientations of $G$ into a
distributive lattice. 

A finite distributive lattice is the lattice of down-sets (also known
as ideals) of some poset $P$. Let a `step' consist in adding/removing
a random element of $P$ to/from the down-set. These step yield the
\emph{up-down Markov chain} on the distributive lattice. A nice
feature of the up-down Markov chain is that it is monotone,
see~\cite{Propp-DL-97}.  A monotone Markov chain is suited for using
\emph{coupling from the past}, see~\cite{PW-exact-96}. This method
allows to sample exactly from the uniform distribution on the elements
of a distributive lattice.

The challenge in applications of the up-down Markov chain is to
analyze its mixing time. In~\cite{Propp-DL-97} some examples of
distributive lattices are described where this chain is rapidly mixing
but there are examples where the mixing is slow.  Miracle, Randall,
Streib and Tetali~\cite{MRST-2016} have investigated the mixing time
of the up-down Markov chain for 3-orientations, a class of
$\alpha$-orientations intimately related to Schnyder woods. They show
that there is a class of plane triangulations such that the up-down Markov chain
on the 3-orientations of these triangulations is slowly mixing. For positive they show
that the chain is rapidly mixing on 3-orientations of 
plane triangulations with maximum degree at most 6.
    
In this paper we present similar results for the up-down Markov chain on
the 2-orientations of plane quadrangulations. These special
2-orientations are of interest because they are related to separating
decompositions, a structure with many applications in floor-planning
and graph drawing. For literature on the subject we refer
to~\cite{dFOdM-tao-01,fhko-blpqr-10,ffno-bsspbc-10} and references given
there. 
Specifically we show that there is a class of plane quadrangulations
such that the up-down Markov chain on
the 2-orientations of these quadrangulations is slowly mixing. 
On the other hand the chain is rapidly mixing on 2-orientations of 
quadrangulations with maximum degree at most 4.

Regarding examples for slow mixing we also revisit the case of
3-orientations, here we have somewhat simpler examples, compared to those
from~\cite{MRST-2016}. Our examples also have a smaller maximum
degree,  $O(\sqrt{n})$ instead of $O(n)$ on graphs with $n$ vertices.
We also exhibit a function~$\alpha$ and a
class of plane graphs of maximum degree 6 such that the up-down Markov
chain on the $\alpha$-orientations of these graphs is slowly mixing.

\section{Preliminaries}

In the first part of this section we give some background on 
the up-down Markov chain on general $\alpha$-orientations.
Then we discuss 2-orientations and the associated separating
decompositions. Finally we provide 
some background on mixing times for Markov chains.

\subsection{The up-down Markov chain of  $\alpha$-orientations}
\label{ssec:up-down}

Let $G$ be a plane graph and $\alpha:V\to \NN$ be such that $G$ admits
$\alpha$-orientations. For $\alpha$-orientations $\vec{G},\vec{G}'$ of
$G$ we define $\vec{G} < \vec{G}'$ whenever $\vec{G}'$ is obtained by
reverting a simple clockwise cycle of $\vec{G}$. This order relation
makes the set of $\alpha$-orientations of $G$ into a distributive
lattice, see~\cite{Fe-04} or~\cite{fk-uld-09}.

The steps of the up-down Markov chain on a distributive lattice 
${\cal L} = (X,<)$
correspond to changes $x \leftrightarrow x'$ for
covering pairs $x \prec x'$, i.e., pairs $x <
x'$ such that there is no $y \in X$ with 
$x < y < x'$. In other words the 
up-down Markov chain performs a random walk on the 
diagram of the lattice.
The transition probabilities are (usually) chosen 
uniformly with a nonzero probability for staying in a state.
Since the diagram of a lattice is connected the chain is ergodic.
It is also symmetric, hence, the unique stationary distribution is
the uniform distribution on the set of $\alpha$-orientations.

The steps of the up-down Markov chain of $\alpha$-orientations are
given by certain reversals of cycles. For a clean description we  need 
the notion of a \term{rigid edge}.  An edge of $G=(V,E)$ is {\em
  $\alpha$-rigid} if it has the same direction in every
$\alpha$-orientation of $G$. Let $R\subseteq E$ be the set of
$\alpha$-rigid edges. Since directed cycles of an $\alpha$-orientation
$\vec{G}$ can be reversed, rigid edges never belong to directed
cycles. Define $r(v)$ as the number of rigid edges that have $v$ as a
tail and let $\alpha'(v) = \alpha(v) - r(v)$.  Now
$\alpha$-orientations of $G$ and $\alpha'$-orientation of $G'=(V,E-R)$
are in bijection. And with the inherited plane embedding of $G'$ the
distributive lattices are isomorphic.  

If $G'$ is disconnected then we can shift connected components of $G'$
to get a plane drawing $G^\#$ without nested components. Since the
orientation, clockwise or counterclockwise, of a directed cycle in
$G'$ and~$G^\#$ is identical the distributive lattices of
$\alpha'$-orientations are isomorphic. The steps of the up-down Markov
chain of $\alpha'$-orientations of~$G^\#$ are easy to describe, they
correspond to the reversal of cycles that form the boundary of bounded
faces, the face boundaries of $G^\#$ are the \term{essential cycles}
for the up-down Markov chain of $\alpha$-orientations of $G$. In
slight abuse of notation we also refer to the up-down Markov chain of
$\alpha$-orientations of $G$ as the \term{face flip Markov chain},
after all the essential cycles of $G$ are faces in~$G^\#$.

A more algebraic description of the lattice for a disconnected $G$ is
as follows: Let $H$ be a component of $G$, then $\cL_\alpha(G)$ 
can  be obtained as the product of
lattices $\cL_{\alpha_1}(G-H)\times \cL_{\alpha_2}(H)$,
where $\alpha_1$ and $\alpha_2$ are the restrictions of $\alpha$
to the vertex sets of $G-H$ and $H$ respectively.

From the previous description it follows that 
the elements of the poset $P_\alpha$ whose down-sets 
correspond to elements of $\cL_\alpha(G)$, i.e., to 
$\alpha$-orientations of $G$, are essential cycles.
It is important to keep the following in mind:

\Fact A. An essential cycle can correspond to several 
elements of the poset $P_\alpha$.
\smallskip

This fact is best illustrated with an
example. Figure~\ref{fig:eulerian}(left) shows the 
octahedron graph $\Goct$ with an Eulerian orientation, this is an $\alpha$
orientation with $\alpha(v)=2$ for all $v$. The orientation is
the minimal one in the lattice, it has no counterclockwise oriented
cycle. Figure~\ref{fig:eulerian}(middle) depicts the poset $P_\alpha$
the labels of the elements of $P_\alpha$ refer to the corresponding
faces of $\Goct$. The elements $1,1',1''$ all correspond to the
same face of $\Goct$, this face has to be reversed three times 
in a sequence of face flips that transforms the minimal
Eulerian orientation into the maximal.  

   \calc_figscale{73}
    \begin{figure}[htb]
    \centerline{\input{\path/eulerian.pstex_t}}
    \caption{\label{fig:eulerian}}
    \end{figure}
    VC
{Left: A minimal $\alpha$-orientation. Middle: The poset
  $P_\alpha$. Right: The $\alpha$-orientation corresponding to the
  down set $\{1,2,3,4,1',6,7,4'\}$ of $P_\alpha$.}

The elements of $P_\alpha$ can be found as follows. 
Let $\Gmin$ be the minimal $\alpha$-orientation, i.e., the one without 
counterclockwise cycles. Starting from $\Gmin$ perform \term{flip}{\em s},
i.e., reversals of essential cycles from clockwise to
counterclockwise, in any order until no further flip is possible.
The unique $\alpha$-orientation that admits no flip is the maximal
one. The flips of a maximal flip-sequence $S$ are the
elements of $P_\alpha$. Let~$\hat{p}(f)$ be the number of times an
essential cycle $f$ has been flipped in $S$. Hence,
the elements of $P_\alpha$ are $\{(f,i) : f$ essential cycle, $1\leq i
\leq \hat{p}(f) \}$. 

If essential cycles $f$ and $f'$ share an edge $e$ then from observing
the orientation of $e$ we find that between any two appearances of $f$
in a flip-sequence there is a appearance of $f'$. From this we obtain

\Fact B. If essential cycles $f$ and $f'$ share an edge, then
$|\hat{p}(f) - \hat{p}(f')| \leq 1$.
\smallskip

The above discussion is based on~\cite{Fe-04} where $\alpha$-orientations of $G$
have been analyzed via $\alpha$-potentials, an encoding of
down-sets of $P_\alpha$. If $\vG$ is an $\alpha$-orientation,
then we say that an essential cycle $f$ is at \term{potential level}
$i$ in $\vG$ if $(f,i)$ belongs to the down-set $D_{\vG}$ of $P_\alpha$
corresponding to $\vG$ but $(f,i+1)\not\in D_{\vG}$.

\subsection{2-orientations and separating decompositions}

A \term{quadrangulation} is a plane graphs whose faces
are uniformly of degree $4$. Equivalently
quadrangulations are maximal bipartite plane graphs.

Let $Q$ be a quadrangulation, we call the color classes of the bipartition
white and black and name the two black vertices on the outer face $s$
and $t$. A \term{2-orientation} of $Q$ is an orientation of the edges
such that $\outdeg(v) = 2$ for all $v\neq s,t$. Since a
quadrangulation with $n$ vertices has $2n-4$ edges it follows that
$s$ and $t$ are sinks.

A \term{separating decomposition} of $Q$ is an orientation and coloring of the
edges of $Q$ with colors red and blue such that two conditions hold:
\Item{(1)} All edges incident to $s$ are ingoing red and all edges incident to
$t$ are ingoing blue.
\Item{(2)} Every vertex $v\neq s,t$ is incident to a nonempty interval of red
edges and a nonempty interval of blue edges. If~$v$ is white, then, in
clockwise order, the first edge in the interval of a color is outgoing and all
the other edges of the interval are incoming. If~$v$ is black, the outgoing
edge is the last one of its color in clockwise order (see
Figure~\ref{fig:vertex-2-cond}). 

   \calc_figscale{25}
    \begin{figure}[htb]
    \centerline{\input{\path/vertex-2-cond.pstex_t}}
    \caption{\label{fig:vertex-2-cond}}
    \end{figure}
    VC
{Edge orientations and colors at white and black vertices.}

\ni 
Separating decompositions have been studied in \cite{BH-2012}, \cite{dFOdM-tao-01},
\cite{ffno-bsspbc-10}, and~\cite{fhko-blpqr-10}. Relevant to us are
the following two facts:

\Fact 1. In a separating decomposition every vertex $v\neq s,t$ 
has a unique directed red path $v \to s$ and a unique blue path $v \to t$. 
The two paths only intersect at $v$.

\Fact 2. The forget function that associates a 2-orientation with a
separating decomposition is a bijection between the set of
2-orientations and the set of separating decompositions of a
quadrangulation.
\medskip

\ni A proof of these facts can be found in~\cite{dFOdM-tao-01}.

\subsection{Markov chains and mixing times}

We refer to~\cite{LevinPeresWilmer} for basics on Markov chains.  In
applications of Markov chains to sampling and approximate counting it
is critical to determine how quickly a Markov chain $M$ converges to
its stationary distribution $\pi$. Let $M^t(x,y)$ be the probability
that the chain started in $x$ has moved to~$y$ in $t$ steps.  The
\term{total variation distance} at time $t$ is 
$
 \| M^t - \pi \|_{TV} = \max_{x\in \Omega} \frac{1}{2} \sum_y
 |M^t(x,y) - \pi(y)|, 
$
here The \term{mixing time} of $M$ is defined as $\Tmix = \min_t (\|
M^t - \pi \|_{TV} \leq 1/4 )$.  The state space $\Omega$ of the Markov
chains considered by us consists of sets of graphs on $n$
vertices. Such a chain is \term{rapidly mixing} if $\Tmix$ is upper
bounded by a polynomial of $n$.

A key tool for lower bounding the mixing time of an ergodic Markov
chain is the \term{conductance} defined as
$ \Phi_M = \min_{S \subseteq \Omega, \pi(S) \leq 1/2}
\frac{1}{\pi(S)} \sum_{s_1 \in S, s_2 \notin  S} \pi(s_1) \cdot
M(s_1,s_2)$.
The connection with $\Tmix$ is given~by
\Fact T.
$\Tmix \geq (4\Phi_M) ^{-1}$.
\medskip

\ni
This is Theorem 7.3 from~\cite{LevinPeresWilmer}. A similar result was already
shown in~\cite{SJ-89}. We will use this inequality mainly in the context 
of \term{hour glass} shaped state spaces where we 
have a partition $\Omega^-,\Omega^0,\Omega^+$ of the state space
with the property that all paths of the transition graph
of the Markov chain that connect $\Omega^-$
and $\Omega^+$ contain a vertex from $\Omega^0$. The following lemma
shows that if  $\Omega^-$
and $\Omega^+$ are large and $\Omega^0$ is small with respect to
$\pi$, then the conductance is small.

\begin{lemma}\label{lem:hour-glass}
If $\Omega^-,\Omega^0,\Omega^+$ is a partition of $\Omega$ such that 
$M(s_1,s_2) = 0$ for all $s_1 \in \Omega^-$ and $s_2 \in \Omega^+$,
then 
$$ 
\Phi_M \leq \frac{\pi(\Omega^0)}{\min\{\pi(\Omega^-),\pi(\Omega^+)\}}.
$$
\end{lemma}
\Proof
We assume that $\pi( \Omega^- ) \leq \pi(\Omega^+)$ and hence
$\pi(\Omega^-) \leq \frac{1}{2}$.
Now 
\begin{eqnarray*}
 \Phi_M & = & \kern-4mm \min_{S \subseteq \Omega, \pi(S) \leq \tfrac{1}{2}} 
        \frac{1}{\pi(S)} \sum_{s_1 \in S, s_2 \notin  S} \pi(s_1) \cdot M(s_1,s_2)
  \quad\leq\quad \frac{1}{\pi(\Omega^-)} \sum_{s_1 \in \Omega^-,
            s_2 \not\in \Omega^-}  \pi(s_1) \cdot M(s_1,s_2) \\
   &=&  \frac{1}{\pi(\Omega^-)} \left(\sum_{s_1 \in \Omega^-, s_2
       \in \Omega^+ } \pi(s_1) \cdot \underbrace{M(s_1,s_2)}_{=0}
     + \sum_{s_1 \in \Omega^-, s_2 \in \Omega^0} \pi(s_1) \cdot
     M(s_1,s_2)\right)\\
   &\leq&  \frac{1}{\pi(\Omega^-)} \left(\sum_{s_2\in \Omega^0 } \sum_{s \in \Omega} \pi(s) \cdot
     M(s,s_2)\right) \quad=\quad  \frac{1}{\pi(\Omega^-)} \sum_{s_2
     \in \Omega^0} \pi(s_2) \\
   &=&
    \frac{\pi(\Omega^0)}{\pi(\Omega^-)} \quad=\quad
    \frac{\pi(\Omega^0)}{\min\{\pi(\Omega^-), \pi(\Omega^+) \}} 
\end{eqnarray*}
\vskip-8mm
\qed

\section{Markov chains for 2-Orientations}

In this section we study the Markov chain $M_2$ for 2-orientations
of plane quadrangulations. This is a special instance of the 
up-down Markov chain for $\alpha$-orientations. A step of the 
chain consists in the reversal of a directed essential cycle.

\begin{lemma}\label{lem:essential4}
The essential cycles for the Markov chain $M_2$ of a plane quadrangulation
are the four-cycles that contain no rigid edge.
\end{lemma}

\Proof
Let $C$ be a four cycle with nonempty interior. We claim that all the edges 
incident to a vertex of~$C$ and a vertex from the interior of $C$ 
are rigid. Let $U$ be the set of vertices interior to $C$ and $E[U]$
be the set of edges incident to a vertex of $U$. Since the cycle
together with $U$ induces a quadrangulation we have 
$|E[U] \cup E_C| = 2|U\cup C| - 4$, i.e., $|E[U]| = 2|U|$.
Hence all edges connecting $U$ to $C$ are out-edges of $U$, this
is the claim.

It follows that every four cycle that contains no rigid edge 
is a face boundary of a component of the non-rigid edges.
This shows that such four-cycles are essential.

Now let $C$ be a cycle of length more than 4 which is a directed cycle in some
2-orientation $\vec{Q}$. A simple counting argument as above shows that in $\vec{Q}$
there is an edge $\vec{e}$ oriented from $C$ into the interior.  From the
correspondence between 2-orientations and separating decompositions together
with Fact~1 we know that there is a directed path $p$ starting with $\vec{e}$
and again hitting $C$. The path~$p$ together with a section of the directed
cycle $C$ is a directed cycle of $\vec{Q}$.  Hence, the edges of $p$ are not
rigid and $C$ is not a boundary of a face of a component of the non-rigid
edges.  \qed

The \term{Markov chain $M_2$} can now be readily described. In each step it
chooses a four-cycle $C$ and $p\in[0,1]$ uniformly at random. If
$C$ is directed in the current orientation $\vec{Q}$ and $p\leq 1/2$, then $C$
is reversed, otherwise the new state equals the old one. The stationary
distribution of $M_2$ is the uniform distribution. (The role of $p$
and the threshold $1/2$ is only to ensure that the Markov chain is
aperiodic.) 

Fehrenbach and R\"uschendorf \cite{FR-2004} have shown that 
$M_2$ is rapidly mixing for certain subsets of the
quadrangular grid. In Subsection~\ref{ssec:rap-mix-2orient} we generalize
this result and prove rapid mixing for quadrangulations 
of maximum degree $\leq 4$.

First, however, we show an exponential lower bound for the mixing time
of $M_2$ on a certain family of
quadrangulations.

\subsection{Slow mixing for 2-orientations}

\begin{figure}[htb]
    \centerline{\input{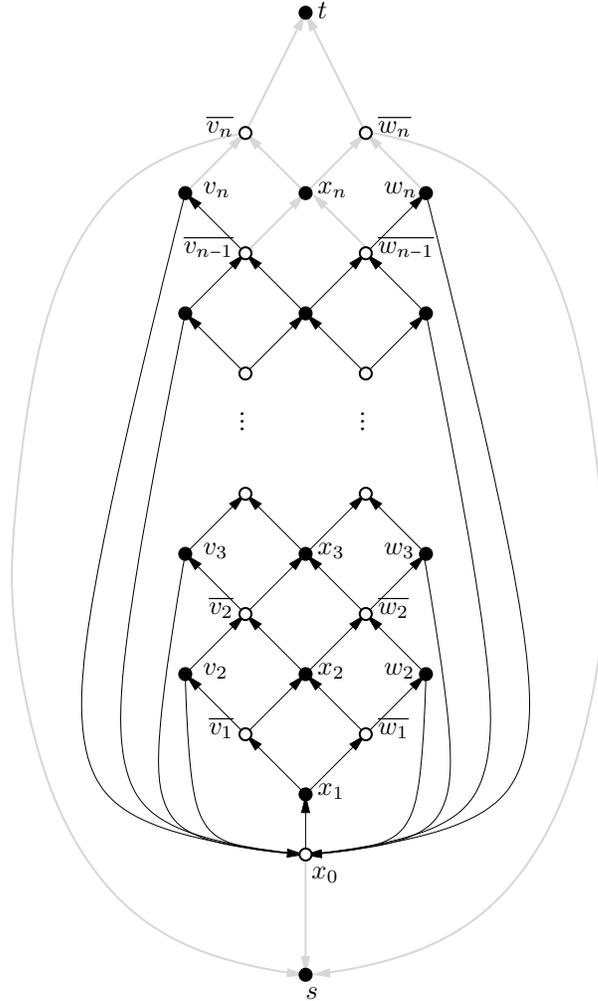}}
    \caption{The graph $Q_n$ with the unique $2$-orientation
  containing the edge $(x_0,x_1)$. Rigid edges
  are shown gray.\label{fig:Qn}}
    \end{figure}

\begin{theorem}\label{thm:lower-2orient}
Let $Q_n$ be the quadrangulation on $5n+1$ vertices shown in
Figure~\ref{fig:Qn}. The  Markov chain $M_2$ on
2-orientations of $Q_n$ has  $\Tmix >  3^{n-3}$.
\end{theorem}
   \calc_figscale{80}
\Proof
Let $\Omega$ be the set of 2-orientations of $Q_n$. We define 
a hour glass partition $\Omega_L,\Omega_c,\Omega_R$ of this
set. The edge $(x_0,s)$ is rigid, the second out-edge $(x_0,a)$ of $x_0$
is called {\em left} if $a\in\{v_2,\ldots,v_n\}$, it is {\em right} if
$a\in\{w_2,\ldots,w_n\}$ and it is {\em central} if $a=x_1$. 
Now $\Omega_L,\Omega_c,\Omega_R$ are the sets 2-orientations where
the second out-edge of $x_0$ is left, central, and right respectively.
With the next claim we show that this is a hour glass partition.

\Claim 1. If $\vec{Q_1}\in \Omega_L$ and $\vec{Q_2} \in \Omega_R$,
then $M_2(\vec{Q_1},\vec{Q_2}) = 0$.
\medskip

\ni
If $\vec{Q} \to \vec{Q}'$ is a step of $M_2$ which changes the 
second out-edge $\vec{e}$ of $x_0$, then the step corresponds to the
reversal of an essential four-cycle containing $\vec{e}$.
Any four-cycle of $Q_n$ that contains $x_0$ either only contains
edges from $x_0$ to vertices from $\{x_1,v_2,\ldots,v_n\}$ or it only
contains edges from $x_0$ to vertices from $\{x_1,w_2,\ldots,w_n\}$.
Hence, if $\vec{Q}\in \Omega_L$, then $\vec{Q}'\in \Omega_L\cup\Omega_c$.
\qedclaim

\Claim 2. $|\Omega_c| = 1$ and Figure~\ref{fig:Qn} shows the unique
2-orientation in this set.
\medskip

\ni
Consider $\vec{Q} \in \Omega_c$.  All the edges between
$\{v_n,x_n,w_n\}$ and $\{\overline{v_n},\overline{w_n}\}$ are oriented
upward in $\vec{Q}$, they are rigid. Suppose all the edges between
$\{v_k,x_k,w_k\}$ and $\{\overline{v_k},\overline{w_k}\}$ are oriented
upward in $\vec{Q}$.  We also know the directed edges $(v_k,x_0)$ and
$(w_k,x_0)$ in $\vec{Q}$. Together this accounts for all out-edges of
$v_k$, $x_k$, and $w_k$. Hence all the edges between
$\{\overline{v_{k-1}},\overline{w_{k-1}}\}$ and $\{v_k,x_k,w_k\}$ are
oriented upward in $\vec{Q}$.  These edges cover all the out-edges of
$\overline{v_{k-1}}$ and $\overline{w_{k-1}}$ whence all edges between
$\{v_{k-1},x_{k-1},w_{k-1}\}$ and
$\{\overline{v_{k-1}},\overline{w_{k-1}}\}$ are oriented upward in
$\vec{Q}$.  With downward induction on $k$ this shows that $\vec{Q}$
has to be the 2-orientation shown in Figure~\ref{fig:Qn}.
\qedclaim

\Claim 3. $|\Omega_L| = |\Omega_R| \geq \frac{1}{2}(3^{n-1}-1)$.
\medskip

\ni
From the symmetry of $Q_n$ we easily get that $|\Omega_L| =
|\Omega_R|$. Now let $P_k$ be the set of directed path from $x_0$ to
$v_k$ in $\vec{Q}$ from Figure~\ref{fig:Qn}. If $p\in P_k$ then
$(v_k,x_0)$ together with $p$ forms a directed cycle in~$\vec{Q}$.
Reverting this cycle yields a 2-orientation that contains the edge
$(x_0,v_k)$. This 2-orientation belongs to $\Omega_L$. Different
paths in $P_k$ yield different orientations. 
Therefore, $|\Omega_L| \geq \sum_k |P_k|$ (in fact equality holds).

It remains to evaluate $|P_k|$. With induction we easily obtain
that in $\vec{Q}$ there are exactly $3^{i-1}$ directed paths from $x_0$ to either
of $\overline{v_i}$ and $\overline{w_i}$. Hence $|P_k| = 3^{k-2}$ and
$|\Omega_L| \geq \sum_{2\leq k\leq n} 3^{k-2} =
\frac{1}{2}(3^{n-1}-1)$.
\qedclaim

The three claims together with Lemma~\ref{lem:hour-glass} yield
$\Phi_{M_2(Q_n)} \leq \frac{2}{3^{n-1}-1}$. Which implies the theorem
via Fact~T.
\qed

\subsection{The tower chain for low degree quadrangulations}
\label{ssec:rap-mix-2orient}

Following ideas originating from~\cite{LuRaSi95} we define a 
tower Markov chain $M_{2T}$ that extends $M_2$. A~single step of 
$M_{2T}$ can combine several steps of $M_2$. Using a coupling
argument we show that $M_{2T}$ is rapidly mixing on 
quadrangulations of degree at most 4. With the comparison
technique this positive result will then be extended to $M_2$.

The basic approach for our analysis of $M_{2T}$ on low degree
quadrangulations is similar to what Fehrenbach and
R\"uschendorf~\cite{FR-2004} did on a class of subgraphs of the
quadrangular grid. In the context of 3-orientations of
triangulations similar methods were applied by
Creed~\cite{Cr-09} to certain subgraphs of the
triangular grid and later by Miracle et al.~\cite{MRST-2016} to
general triangulations.
As Creed~\cite{Cr-09} noted there is an inaccurate claim in the proof 
of~\cite{FR-2004}. Later Miracle et al.~\cite{MRST-2011} stepped into the same trap
(it has been corrected in the final version~\cite{MRST-2016}). 
In~\ref{sssec:problems} below we discuss 
these issues and show how to repair the proofs. 

Let $Q$ be a quadrangulation on $n$ vertices and $\Omega$ be
the set of 2-orientations of $Q$. From the considerations in
Subsection~\ref{ssec:up-down} we know that there is 
a redrawing $Q^\#$ of the subgraph of non-rigid edges of $Q$
such that the essential cycles of $Q$ are the boundaries 
of bounded faces of $Q^\#$. From
Lemma~\ref{lem:essential4} we know that these faces are four-faces.

Let $\vec{Q}$ be a 2-orientation and $C$ be a simple cycle.  With
$e^+(C)$ we denote the number of clockwise edges of $C$ and with
$e^-(C)$ the number of counterclockwise edges.  If $f$ is a four-cycle
and $\nu(f) = |e^+(f) - e^-(f)|$, then $\nu(f)$ can take the values
$0$, $2$, and $4$. The face $f$ is \nct{oriented} if $\nu(f)=4$, it is
\term{scrambled} if $\nu(f)=0$, and it is \term{blocked} is
$\nu(f)=2$.  If $f$ is blocked, then three edges have the same
orientation and one edge does not. We call this the \nct{blocking
  edge} of $f$.

\begin{figure}[htb]
    \centerline{\input{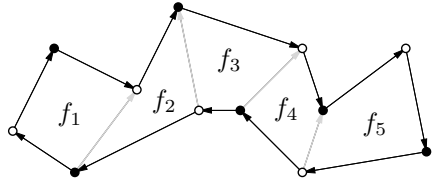}}
    \caption{A tower of length 5.\label{fig:tower}}
    \end{figure}

A \term{tower} of length $k$ is a sequence $(f_1,f_2,\ldots,f_k)$ of
four-cycles of $\vec{Q}$ such that each $f_i$ for $i=1,..,k-1$ is
blocked and $f_k$ is oriented. Moreover, in $f_i$ the blocking edge of
$f_{i-1}$ is opposite to the blocking edge of $f_i$ for $i=1,..,k-1$. A
tower of length $1$ is just an oriented face, Figure~\ref{fig:tower}
shows a tower of length 5.

   \calc_figscale{55}

\vbox{}\vskip-10mm\vbox{}
Lemma~\ref{lem:simple-tower} below implies that that removing all
blocking edges from a tower $T$ of length $k$ we obtain a connected
region whose boundary $\partial T$ is an oriented cycle with $2k+2$
edges. This is the \term{boundary cycle} of the tower. The boundary
cycle need not be simple but each edge of $\partial T$ only belongs to
a single face $f_i\in T$. Therefore, we can also obtain the effect of
reverting $\partial T$ by reverting $f_k,f_{k-1},\ldots,f_1$ in this
order.

For the following arguments we assume that $Q$ has no nested
four-cycles. This is justified by the isomorphism between the lattices 
of 2-orientations of $Q$ and of $\alpha'$-orientations of $Q^\#$.
\vskip-6mm\vbox{}
\begin{lemma}\label{lem:simple-tower}
If $f_i$ and $f_j$ are different faces of a tower $T$ and they share an edge,
then $j\in \{i-1,i+1\}$ and the shared edge is the blocking edge of
one of them. 
\end{lemma}
\Proof The construction sequence $f_1,\ldots,f_k$ of a tower
$T=(f_1,\ldots,f_k)$ yields a directed walk in the dual. The edges of this
walk are the duals of the blocking edges. Each $f$ has at most one blocking
edge and $f_k$ has no blocking edge. Hence, there is no repetition of faces in
a tower. It follows that $\partial T$ is an oriented cycle. Two faces $f_i\neq
f_j$ do not share an edge $e$ of $\partial T$. This is because they would be
the faces of the two sides of $e$ whence $e$ would be clockwise in one of them
and counterclockwise in the other, however, $\partial T$ is uniformly
oriented.\qed
\vskip-6mm\vbox{}
\begin{lemma}\label{lem:unique-tower}
If $f$ is a four-cycle, then there is at most one tower starting with
$f=f_1$. 
\end{lemma}
\Proof
Again we look at the construction sequence of a tower and the 
corresponding directed path in the dual. Each $f_i$ has at most one
blocking edge, hence, there is a unique candidate for $f_{i+1}$.
If~$f_{i+1}$ is oriented it completes the tower. If~$f_{i+1}$
is blocked and the blocking edge is opposite to the edge shared with
$f_i$ the construction of the potential tower can be extended.
Otherwise, there is no tower starting with $f$.
\qed

We are ready to describe the \term{tower Markov chain $M_{2T}$}.
If $M_{2T}$ is in state $\vX$ then it performs the transition to the next step
as follows: an essential four-cycle~$f$, and a $p\in[0,1]$ are each chosen
uniformly at random.  If in $\vX$ there is a tower $T_f$ of length $k$
starting with $f$ then revert $\partial T_f$ if 

\Bitem $\partial T_f$ is clockwise and either $k=1$ and $p\leq 1/2$ or $k>1$
and $p \leq 1/(4k)$, 

\Bitem $\partial T_f$ is counterclockwise and either $k=1$ and $p < 1/2$ or
$k>1$ and $p \geq 1-1/(4k)$. 
\medskip

\ni
In all other cases the new state is again $\vX$.
\bigskip\bigbreak

Since the steps of $M_2$ are also steps of $M_{2T}$ 
the chain is connected. In the orientation 
obtained by reverting the tower 
$T=(f_1,\ldots,f_k)$ there is the tower 
$T'=(f_k,\ldots,f_1)$ whose reversal leads back to the original
orientation. Since both towers have the same length
the chain is symmetric and its stationary distribution is uniform.

The next lemma is where the degree condition is indispensable.

\begin{lemma}\label{lem:clean-tower}
  Let $Q$ have maximum degree $\leq 4$ and let $T=(f_1,\ldots,f_k)$
  be a tower and $\hf\neq f_k$ be an oriented face in a
  $2$-orientation $\vec{Q}$
  of $Q$. If $T$ and $\hf$ share an edge $e$ but $\hf$ and $f_1$
  share no edge, then $e$ is the edge of $f_k$ opposite to the
  blocking edge of $f_{k-1}$.
\end{lemma}

\Proof Let $(u_i,v_i)$ be the blocking edge of $f_i$. We extend the
labeling of vertices of $T$ such that~$\partial T$ is the directed cycle
$v_0,v_1,\ldots,v_{k-1},v_k,u_k,u_{k-1},\ldots,u_1,u_0$.

If $(u_{i+1},u_i)$ with $i\geq 1$ is an edge of $\hf$ and $u_{i-1}\not\in
\hf$, then $\hf$ contains an out-edge of $u_i$ which is not part of $T$.
However, $u_i$ contains the out-edges $(u_i,v_i)$ and $(u_i,u_{i-1})$.  This
contradicts $\outdeg(u_i) = 2$.

If $(v_i,v_{i+1})$ with $i\geq 1$ is an edge of $\hf$ and $v_{i-1}\not\in
\hf$, then $\hf$ contains an in-edge of $v_i$ which is not part of $T$.
Vertex $v_i$ also contains the in-edges $(u_i,v_i)$ and $(v_{i-1},v_i)$. Now
$v_i$ has in-degree~$\geq 3$, since $\outdeg(v_i) = 2$ the degree is at
least~5. A contradiction.

We are not interested in edges shared by $\hf$ and $f_1$, i.e, in edges
containing $u_0$ or $v_0$. Therefore, the only remaining candidate for $e$ is
the edge $(v_k,u_k)$.  \qed

\begin{theorem}\label{thm:mix-2-orient}
Let $Q$ be a plane quadrangulation with $n$ vertices so that each inner
vertex is adjacent to at most 4 edges. The mixing time of
$M_{2T}$ on 2-orientations of $Q$ satisfies $\Tmix\in O(n^5)$.
\end{theorem}

The proof of Theorem~\ref{thm:mix-2-orient} is based on the path coupling 
theorem of Dyer and Greenhill~\cite{DG-1998}. Before stating a simple
version of the Dyer--Greenhill Theorem we need a definition.
A \emph{coupling} for a Markov chain $M$ on a state space $\Omega$ is
a pair $(X_t,Y_t)$ of processes satisfying two conditions:
\Bitem Each of $(X_t)$ and $(Y_t)$ represents $M$, i.e., 
$\Pr(Z_{t+1} = j | Z_t=i) = M_{i,j}$, for $Z\in \{X,Y\}$ and all $t$. 
\Bitem If $X_t=Y_t$ then $X_{t+1} = Y_{t+1}$.  

\vbox{
\begin{theorem}[Dyer--Greenhill]\label{thm:DGcoupling}
Let  $M$ be a Markov chain with state space
$\Omega$. If there is a graph ${\cal G}_M$ with vertex set $\Omega$ 
and a coupling $(X_t,Y_t)$ of $M$ such that with the graph distance 
$d : \Omega \times \Omega \rightarrow \NN$ based on ${\cal G}_M$ we have:
\smallskip

\centerline{$ 
\E[d(X_{t+1},Y_{t+1})] \leq d(X_t,Y_t) \text{\quad and\quad} 
\Pr( d(X_{t+1},Y_{t+1}) \neq d(X_{t},Y_{t})) \geq \rho $}
\smallskip

\ni
then $\Tmix(M) \leq 2\lceil{e}/{\rho}\rceil 
\text{diam}({\cal G}_M)^2$.
\end{theorem}
}

The coupling of $M_{2T}$ used for the proof of Theorem~\ref{thm:mix-2-orient}
is the trivial one, i.e., we run chains $X_t$ and $Y_t$ with the same choices
of $f$ and $p$ in each step.

The graph ${\cal G}$ will be the transition graph of $M_2$, i.e, the
distance between 2-orientations $\vec{X}$ and~$\vec{Y}$ equals the
number of four-cycles that have to be reversed to get from $\vec{X}$
to~$\vec{Y}$.

\begin{lemma}\label{lem:max-pot}
The maximum potential $\hat{p}_\text{max} = \max_f \hat{p}(f)$ of an essential
cycle is less than $n$.
\end{lemma}
\Proof 
Let $Q$ be the quadrangulation whose 2-orientations are in question.
It is convenient to replace $Q$ by $Q^\#$ so that essential cycles are
just faces.
Recall that $\hat{p}$ of the outer face is
$0$ and $|\hat{p}(f) - \hat{p}(f')| \leq 1$ for any two
adjacent faces (Fact~B). Since a quadrangulation has $n-2$ faces we obtain 
$(n-3)$ as an upper bound for $\hat{p}_\text{max}$.
\qed

\begin{lemma}\label{lem:dist(G)}
The diameter of ${\cal G}$ is at most $n^2/2$.
\end{lemma}
\Proof 
The height of the lattice  ${\cal L}_\alpha(Q^\#)$ is the length of a
maximal flip sequence, i.e., 
$\sum_f \hat{p}(f)$. Using (Fact~B) as in the proof of the previous
lemma we find that $\sum_f \hat{p}(f) \leq 0+1+\ldots+(n-3)$. 
This is $< n^2/2$.

In the diagram of a distributive lattice the diameter is attained by the
distance between the zero and the one, i.e., between the global minimum and
the global maximum.  This distance is exactly the height of the lattice.
Since $\cal G$ is the cover graph (undirected diagram) of the distributive
lattice~${\cal L}_\alpha(Q)$ we obtain that the diameter of $\cal G$ is at
most $n^2/2$.  \qed

\subsubsection{Finding an appropriate $\rho$}
\label{sssec:problems}

\begin{wrapfigure}[16]{r}{0.33\textwidth}
\vskip-6mm
\centering
\includegraphics[width=0.30\textwidth]{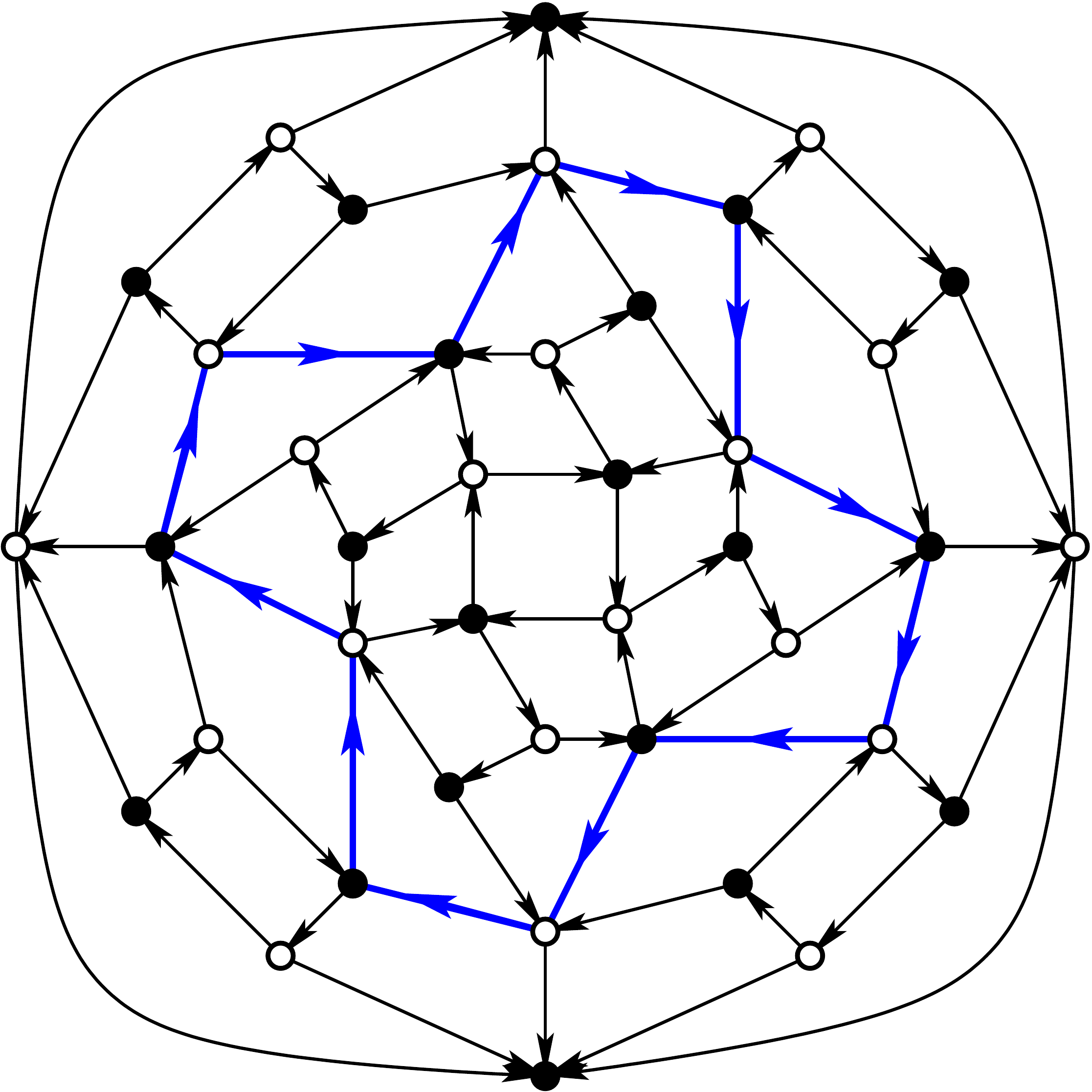}
\vskip-2mm
\caption{\label{fig:hardness-ex}
A quadrangulation $Q$ with a 2-orientation $\vec{Q}$. The 2-orientation
$\vec{Q}'$ obtained by reverting the blue cycle has the same oriented
faces. Hence, all face reversals preserve the distance.}
\end{wrapfigure}
To get a reasonable $\rho$ the following argument is tempting and was
actually used in~\cite{FR-2004} and~\cite{MRST-2011}: For given $\vX$
and $\vY$ there is always at least one essential cycle~$f$ whose reversal
in $\vec{X}$ reduces the distance to $\vec{Y}$.  If $(X_t,Y_t) =
(X,Y)$ and this cycle $f$ is chosen by~$M_{2T}$, then with probability
$1/2$ the distance is reduced. There are at most $n-3$
essential cycles.  Hence we may set $\rho = 1/(2n)$.


Indeed for up-down Markov chains on distributive such a statement
holds. If $I$ and $J$ are down-sets of the poset $P$, then there is 
an $x\in P$ whose addition to or removal from $I$ decreases the
distance to $J$.
In the context of $\alpha$-orientations, however, an $f$ whose reversal in
$\vec{X}$ reduces the distance to~$\vec{Y}$ may be oriented in~$\vY$ with the
same orientation as in $\vX$.  In that case if $f$ is chosen by $M_{2T}$ the
reversal of $f$ is applied to both or to none Figure~\ref{fig:hardness-ex}
shows that there are cases where a pairs $(X_t,Y_t)$ exists such that 
$\Pr(d(X_{t+1},Y_{t+1}) \neq d(X_{t},Y_{t})) = 0$.

To overcome this problem we now define the \term{slow tower Markov
  chain} $M_{S2T}$. 
\smallskip

\ni
If $M_{S2T}$ is in state $\vX$ then it performs the transition to the
next step as follows:
an essential four-cycle $f$, a value $i$ with $0\leq i < n$
 and a $p\in[0,1]$ are each chosen uniformly at random. 
If $f$ is not at potential level~$i$ in $\vX$, then nothing is done
and $\vX$ is the new state.

Otherwise, if there is a tower $T_f$ of length $k$ starting with $f$ then revert $\partial T_f$
if 
\Bitem $\partial T_f$ is clockwise and either $k=1$ and $p\leq 1/2$ or
$k>1$ and $p \leq 1/(4k)$,
\Bitem $\partial T_f$ is counterclockwise and either $k=1$ and $p > 1/2$ or
$k>1$ and $p \geq 1-1/(4k)$.
\smallskip

In all other cases the new state is again $\vX$.
\bigskip

\begin{lemma}\label{lem:slow-rho}
If $(X_t,Y_t)$ is a trivial coupling of the slow chain $M_{S2T}$, then
$$\Pr( d(X_{t+1},Y_{t+1}) \neq d(X_{t},Y_{t})) \geq 1/(2n^2).$$
\end{lemma}

\Proof
For given $\vX$ and $\vY$ there is always at least one essential cycle~$f_1$ whose reversal
in $\vec{X}$ reduces the distance to $\vec{Y}$.  
If $f_1$ appears in $\vX$ and $\vY$ with the same orientation 
then the potential level of $f_1$ in $\vX$ and~$\vY$ is different.
Hence, if for the step of $M_{S2T}$ the triple $(f,i,p)$ is 
chosen such that $f=f_1$ and $i$ is the potential level of $f$ in
$\vX$ and $p$ is such that $f$ is actually reversed, then the 
distance decreases.

The probability for choosing $f$ is at least $1/n$. For $i$ and~$p$ 
the probabilities are $1/n$ and $1/2$ respectively.
Together this yields the claimed bound.
\qed

\subsubsection{Completing the proof of Theorem~\ref{thm:mix-2-orient}}

In Lemma~\ref{lem:change-dist} we show that if $(\vX, \vY)$ is an edge
of ${\cal G}$ and $(\vX^+, \vY^+)$ is the pair obtained after a single
coupled step of the tower chain $M_{2T}$, then $\E[d(\vX^+, \vY^+) -
d(\vX,\vY)] \leq 0$.  Note that a step of the coupled slow
chain~$M_{S2T}$ moves the pair $(\vX,\vY)$ to $(\vX^+,\vY^+)$ with
probability $1/n$ and otherwise stays at $(\vX,\vY)$. Hence
Lemma~\ref{lem:change-dist} also applies to $M_{S2T}$.

Assuming Lemma~\ref{lem:change-dist} we get the following:
\begin{proposition}
Let $Q$ be a plane quadrangulation with $n$ vertices so that each inner
vertex is adjacent to at most 4 edges. The mixing time of
$M_{S2T}$ on 2-orientations of $Q$ satisfies $\Tmix(M_{S2T})\in O(n^6)$
\end{proposition}

\Proof For the condition $ \E[d(X_{t+1},Y_{t+1})] \leq d(X_t,Y_t)$ needed
for the application of Theorem~\ref{thm:DGcoupling} 
we need Lemma~\ref{lem:change-dist}. The inequality from the lemma is
also true for $M_{S2T}$ because this behaves like a slowed down
version of $M_{2T}$. Linearity of expectation allows to transfer the 
inequality from single edges to paths.  

An application of Theorem~\ref{thm:DGcoupling} with parameters $\rho =
\frac{1}{2n^2}$ (Lemma~\ref{lem:slow-rho}) and $\text{diam}({\cal G})
\leq n^2/2$ (Lemma~\ref{lem:dist(G)}) yields $\Tmix(M_{S2T}) \leq
en^6$.  \qed

The mixing time of the slow chain could thus be proven with a coupling
that allows an application of the theorem of Dyer and Greenhill.
Now consider a single state $\vX_t$ evolving according to the
slow chain $M_{S2T}$. Note that this is exactly as if we would run 
the tower chain $M_{2T}$ but only allow a transition to be 
conducted if an additional uniform random variable $q \in
\{0,\ldots,n-1\}$ takes the value $q=0$.
It follows that the mixing times of $M_{S2T}$ and
of $M_{2T}$ deviate by a factor of $n$. Therefore, 
$\Tmix(M_{2T}) \leq en^5$. 

To complete the proof of Theorem~\ref{thm:mix-2-orient} it remains to
prove Lemma~\ref{lem:change-dist}.

\begin{lemma}\label{lem:change-dist}
  If $(\vX, \vY)$ is an edge of ${\cal G}$ and $(\vX^+, \vY^+)$ is
  the pair obtained after a single coupled step of~$M_{2T}$, then
  $\E[d(\vX^+, \vY^+) - d(\vX,\vY)] \leq 0$.
\end{lemma}

\Proof
Since $(\vec{X}, \vec{Y})$ is an edge of ${\cal G}$ they differ
in the orientation of exactly one face $\hf$. We assume w.l.o.g\ that
$\hf$ is oriented clockwise in $\vX$ and counterclockwise in $\vY$.

Let $f$ be the face chosen for the step of $M_{2T}$. Depending on
$f$ we analyze $d(\vX^+,\vY^+)$ in three cases.

\ni{\bfseries A.}
If $f=\hf$, then depending on the value of $p$ face $f$ is
reversed either in $\vX$ or in $\vY$.
After the step the orientations $\vX^+$, $\vY^+$ coincide.
The expected change of distance in this case is $-1$.

\ni{\bfseries B.}
If $f$ and $\hf$ share an edge and $f\neq \hf$ there are three
options depending on the type of $f$ in~$\vY$. 

{\bfseries 1.}  Face $f$ is oriented in $\vY$, necessarily clockwise. It
follows that in $\vX$ face $f$ starts the clockwise tower $(f,\hf)$ of length
two. In $\vY$ a face $f$ is a clockwise tower of length 1. If $p \leq 1/8$
both towers are reversed so that $\vX^+$ and $\vY^+$ coincide. If $1/8 < p
\leq 1/2$, then $f$ is reversed in $\vY$ while $\vX^+=\vX$, in this case the
distance increases by 1. If $p > 1/2$ both orientations remain unchanged.  The
expected change of distance in this case is $\tfrac{1}{8}\cdot ( - 1) +
(\tfrac{1}{2} - \tfrac{1}{8})\cdot(+1) + \tfrac{1}{2} \cdot 0 = \tfrac{1}{4}$.

{\bfseries 2.}  Face $f$ is scrambled in $\vY$. In this case $f$ is
blocked in $\vX$ and it may start a tower of length~$k$. If $p\leq
1/(4k)$ this tower is reverted which results in a increase of distance
by $k$. In all other cases the distance remains unchanged. Hence, the
expected change of distance in this case is $\leq 1/4$.

{\bfseries 3.}  Face $f$ is blocked in $\vY$. Then it is either
oriented or scrambled in $\vX$. After changing the role of $\vX$ and
$\vY$ we can use the analysis of the other two cases to conclude that
the expected change of distance is again $\leq 1/4$.

\ni{\bfseries C.}
Finally, suppose that $f$ and $\hf$ have no edge in common. 

{\bfseries 1.} 
If $f$ starts a tower in $\vX$ which has no edge in common
with $\hf$, then $f$ starts the very same tower in~$\vY$ and the 
coupled chain will either revert both towers or none of them. 
The distance remains unchanged.

{\bfseries 2.}  Now let $f$ start a tower $T=(f_1,\ldots,f_k)$ in
$\vX$ which has an edge in common with~$\hf$. The case where~$\hf$ and
$f_1=f$ share an edge was considered in {\bfseries B}. Now
Lemma~\ref{lem:clean-tower} implies that either $\hf=f_k$ or
$\hf\neq f_k$ and the shared edge is such that $(f_1,\ldots,f_k,\hf)$
is a tower in $\vY$. Hence, with $T$ there is a tower~$T'$ in $\vY$
that starts in $f$ and has length $k\pm 1$, moreover $T$ and $T'$ have
the same orientation. Let $\ell$ be the larger of the lengths of $T$
and $T'$. With a probability of $1/(4\ell)$ both towers are reversed and
the distance decreases by 1.  With a probability of $1/(4(\ell-1)) -
1/(4\ell)$ only the shorter of the two towers is reversed and the
distance increases by $\ell-1$.  With the remaining probability both
orientations remain unchanged.  The expected change of distance in
this case is $\tfrac{1}{4\ell}\cdot ( - 1) + (\tfrac{1}{4(\ell-1)} -
\tfrac{1}{4\ell})\cdot(\ell-1) = 0$.  
\medskip

Let $m$ be the number of essential four-cycles, i.e., the number of
options for $f$. Combining the values for the change of distance in 
cases {\bfseries A, B, C} and the probability of these cases we
obtain:
$$\E[d(\vX^+, \vY^+) - d(\vX,\vY)] \leq \frac{1}{m}(-1) +
\frac{4}{m}(1/4) + \frac{m-5}{m}0 = 0.$$
\vskip-18pt\qed\medskip

\subsection{Comparison of $M_{2T}$ and $M_2$}
\def\AA{{\cal A}}

The comparison of the mixing times of $M_{2T}$ and $M_2$ is based on
a technique developed by Diaconis and Saloff-Coste \cite{DSC-1993}. 
We will use Theorem~\ref{thm:RTcomparison} a variant due to
Randall and Tetali~\cite{RT-1997}.

Let $M$ and $\widetilde{M}$ be two reversible Markov chains on the
same state space $\Omega$ such that $M$ and~$\widetilde{M}$ have the
same stationary distribution $\pi$. With $E(M)$ we denote the edges of
the directed transition graph of~$M$, i.e, $(x,y)\in E(M)$ whenever
$M(x,y) > 0$. Define $E(\widetilde{M})$ alike. 
For each $(x,y)\in E(\widetilde{M})$ define a canonical path
$\gamma_{xy}$ as a sequence $x= v_0, v_1, \dots, v_k = y$ of
transitions of~$M$, i.e. $(v_i, v_{i+1}) \in E(M)$ for all~$i$.
Let $|\gamma_{xy}|$ be the length of $\gamma_{xy}$ and for $(x,y)\in
E(M)$ let
$\Gamma(x,y) := \{ (u,v) \in E(\widetilde{M}) :  (x,y) \in \gamma_{uv} \}.$
Further let
$$ \AA := \max_{(x,y) \in E(M)} {}
       \left\{  \frac{1}{\pi(x)M(x,y)} 
          \sum_{(u,v)\in\Gamma(x,y)}|\gamma_{uv}|\pi(u)\widetilde{M}(u,v)\right\}
$$
and let $\pi_\star := \min_{x \in \Omega} \pi(x)$.

\begin{theorem}[Randall--Tetali]\label{thm:RTcomparison}
In the above setting $\Tmix(M) \leq
4\log(4/\pi_\star)\;\AA\;\Tmix(\widetilde{M})$.
\end{theorem}

We are going to apply this theorem with $M = M_2$ and
$\widetilde{M}=M_{2T}$. Both chains are symmetric, hence reversible,
and have the uniform distribution $\pi$ as stationary distribution.

The definition of the canonical paths comes quite natural.
A transition $(\vec{U},\vec{V})$ of $M_{2T}$ corresponds to the
reversal of  $\partial T$ for some tower $T$ of $\vec{U}$. 
Suppose that $T=(f_1,\ldots,f_k)$ and recall that
the effect of reverting $\partial T$ can also
be obtained by reverting $f_k,f_{k-1},\ldots,f_1$ in this order.
Reverting them one by one yields a path in $E(M)$, this path is chosen to 
be $\gamma_{\vec{U}\vec{V}}$. 

If $|\gamma_{\vec{U}\vec{V}}|=k$, i.e., the transition
$(\vec{U},\vec{V})$ corresponds to a tower of length $k$, then
$M_{2T}(\vec{U},\vec{V}) = 1/(4k)$, hence,
$|\gamma_{\vec{U}\vec{V}}|M_{2T}(\vec{U},\vec{V}) = 1/4$. Also $\pi$
is constant so that $\pi(\vec{U})/\pi(\vX) =1$. For an upper bound on~$\AA$
we therefore only have to estimate the number of tower moves that 
have a canonical path
that contains the face flip at $f$ that moves $\vX$ to $\vY$. If
$T=(f_1,\ldots,f_k)$ is such a tower with $f=f_i$, then
$(f_1,\ldots,f_{i-1},f)$ is a tower in $\vX$ and
$(f_k,\ldots,f_{i+1},f)$ is a tower in $\vY$. Since a tower is defined
by its initial face each of $\vX$ and~$\vY$ has at most~$n$ towers,
all the more each has at most $n$ towers ending in $f$.  This shows
$|\Gamma(\vX,\vY)| \leq n^2$ and $\AA \leq n^2/4$.

It remains to find $\pi_\star=\frac{1}{|\Omega|}$. Since a
quadrangulation has $2n-4$ edges it has at most $2^{2n}$ orientations
this would suffice for our purposes. However, a better upper bound of $1.9^n$
for the number of 2-orientations was obtained in~\cite{FZ-08}.

Given the above ingredients for the comparison theorem 
and the mixing time of $\Tmix(M_{2T})\in O(n^5)$ from
Theorem~\ref{thm:mix-2-orient} we finally have shown
rapid mixing for $M_2$ on certain quadrangulations.

\begin{theorem}
Let $Q$ be a plane quadrangulation with $n$ vertices so that each inner
vertex is adjacent to at most 4 edges. The mixing time of the 
face reversal Markov chain $M_{2}$ on 2-orientations of $Q$ satisfies 
$\Tmix(M_2)\in O(n^8)$.
\end{theorem}

\section{Slow mixing for 3-orientations}

A \term{triangulation} is a plane graphs whose faces
are uniformly of degree $3$. Equivalently
triangulations are maximal plane graphs.

A \term{3-orientation} of a triangulation $T$ is an orientation of the 
internal edges, i.e., of the edges except the three edges of the outer face,
such that $\outdeg(v) = 3$ for all inner vertices $v$. Since a
triangulation with $n$ vertices has $3n-9$ inner edges it follows that
the three outer vertices are sinks. 

A \term{Schnyder wood} of $T$ is an orientation and coloring of the
edges of $T$ with colors red, green, and blue such that two conditions hold:
\Item{(1)} If the vertices of the outer face are colored red, green
and blue in clockwise order, then all inner edges incident to a vertex
$s$ of the outer face are oriented towards $s$ and colored in the
color of $s$.
\Item{(2)} Every inner vertex $v$ has three outgoing edges 
colored red, green, and blue in clockwise order. Incoming 
edges in the sector between two outgoing edges are colored in the
third color (see Figure~\ref{fig:vertex-2-cond}).

   \calc_figscale{25}
    \begin{figure}[htb]
    \centerline{\input{\path/node-prop.pstex_t}}
    \caption{\label{fig:node-prop}}
    \end{figure}
    VC
{The two conditions for Schnyder woods.}

\ni Schnyder woods were introduced by Schnyder in~\cite{s-pgpd-89}.
We refer to~\cite{dFOdM-tao-01,ps-ocst-06,f-gga-04,abfkku-ccoc-13} and
the references given there for properties, applications and
generalizations of Schnyder woods.  Relevant to us is the following
fact, see~\cite{dFOdM-tao-01}:

\Fact 3. The forget function that associates a 3-orientation with a
Schnyder wood is a bijection between the set of 3-orientations and the 
set of Schnyder woods of a triangulation. 
\medskip

From the correspondence between Schnyder woods and $3$-orientations 
it follows that the triangle flip Markov chain can be used to sample 
from either of these structures. The mixing time of this Markov
chain was studied by Creed~\cite{Cr-09} for certain subgraphs of the 
triangular grid and then By Miracle et al.~\cite{MRST-2016} for 
general triangulations. Here we want to revisit the following negative
result. 
 
\begin{theorem}[Miracle--Randall--Streib--Tetali]\label{thm:MRSTlower}
  There is a triangulation $T'_n$ with $4n+1$ vertices with maximum degree
  $2n+3$ such that the triangle flip Markov chain $M_3$ on 3-orientations of
  $T'_n$ has $\Tmix > \frac{1}{16}2^{n/4}$.
\end{theorem}

With Theorem~\ref{thm:lower-3orient} we prove a similar result with a larger
exponential bound on $\Tmix$.  Moreover, $T_n$ is simpler than $T'_n$. This
carries over to the simplicity of the proof. In fact the proof is very similar
to the proof for Theorem~\ref{thm:lower-2orient}.  Below, in
Proposition~\ref{prop:degree} we modify $T_n$ to show that slow mixing of the
triangle flip chain $M_3$ can also be observed for triangulations with maximum
degree in the order of $\sqrt{n}$.

\begin{figure}[htb]
    \centerline{\input{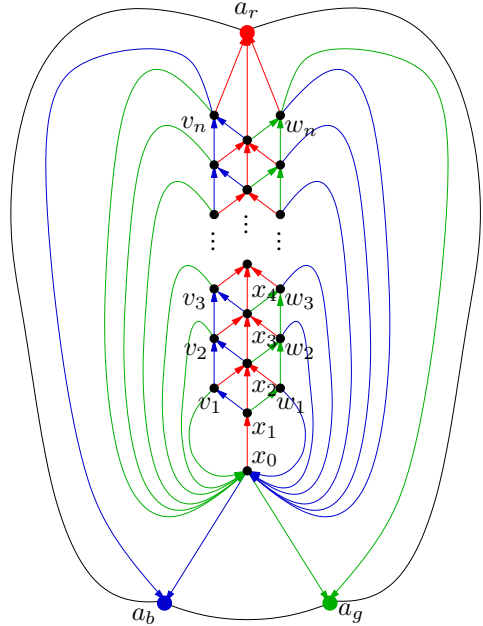}}
    \caption{The triangulation $T_n$ with a Schnyder wood.\label{fig:Tn}}
    \end{figure}

\begin{theorem}\label{thm:lower-3orient}
  Let $T_n$ be the triangulation on $3n+4$ vertices with maximum degree $2n+3$
  shown in Figure~\ref{fig:Tn}.  The triangle flip Markov chain $M_3$ on
  3-orientations of $T_n$ has $\Tmix > (2+\sqrt{3})^{n-2} \approx 3.732^{n-2}$.
\end{theorem}

   \calc_figscale{55}
\Proof
Let $\Omega$ be the set of Schnyder woods of $T_n$. We define 
a hour glass partition $\Omega_L,\Omega_c,\Omega_R$ of this
set. The edges $(x_0,a_g)$ and $(x_0,a_b)$ are rigid, 
the red out-edge $(x_0,z)$ of $x_0$
is called {\em left} if $z\in\{v_1,\ldots,v_n\}$, it is {\em right} if
$z\in\{w_1,\ldots,w_n\}$ and it is {\em central} if $z=x_1$. 
Now $\Omega_L,\Omega_c,\Omega_R$ are the sets Schnyder woods where
the red edge of $x_0$ is left, central, and right respectively.
With the next claim we show that this is a hour glass partition.

\Claim 1. If ${S_1}\in \Omega_L$ and ${S_2} \in \Omega_R$,
the $M_2({S_1},{S_2}) = 0$.
\medskip

\ni
If ${S} \to {S}'$ is a step of $M_3$ which changes the 
red out-edge $\vec{e}$ of $x_0$, then the step corresponds to the
reversal of a triangle containing $\vec{e}$.
There is no triangle in $T_n$ with vertices $\{x_0,v_i,w_j\}$ for
$i,j \in [n]$.
Hence, if ${S}\in \Omega_L$, then ${S}'\in \Omega_L\cup\Omega_c$.
\qedclaim

\Claim 2. $|\Omega_c| = 1$ and Figure~\ref{fig:Tn} shows the 
unique Schnyder wood of this set.
\medskip

\ni
Consider ${S} \in \Omega_c$.  From $(x_0,x_1) \in S$ 
we conclude that $\{v_1,x_2,w_1\}$ are the out-neighbors of $x_1$.
From the degrees it follows that all the edges between $\{v_1,x_2,w_1\}$
and $\{v_2,x_3,w_2\}$ are oriented upward in~$S$. 
Inductively we find that all the edges between $\{v_{i-1},x_i,w_{i-1}\}$
and $\{v_i,x_{i+1},w_i\}$ are oriented upward in~$S$. Since the edges 
$(v_i,x_0)$ and $(w_i,x_0)$ are in $S$ anyway it follows that the 
orientation of all edges is fixed when $(x_0,x_1)$ is fixed.
The bijection between $3$-orientations and Schnyder woods then yields
that the Schnyder wood shown in Figure~\ref{fig:Tn} is the unique element
of $\Omega_c$.
\qedclaim

\Claim 3. $|\Omega_L| = |\Omega_R| > (2+\sqrt{3})^{n-1}$.
\medskip

\ni From the symmetry of $T_n$ we easily get that $|\Omega_L| =
|\Omega_R|$. Now let $P_k$ be the set of directed path from~$x_0$
to~$v_k$ in the orientation $S$ from Figure~\ref{fig:Tn}. If $p\in
P_k$ then $(v_k,x_0)$ together with $p$ forms a directed cycle in~$S$.
Reverting this cycle yields a 3-orientation that contains the edge
$(x_0,v_k)$. This 3-orientation belongs to~$\Omega_L$. Different paths
in $P_k$ yield different orientations.  Therefore, $|\Omega_L| \geq
\sum_k |P_k|$ (in fact equality holds).

It remains to evaluate $g_k = |P_k|$. To do so let $h_k$ be the number of
directed paths from $x_0$ to $x_k$. Clearly, $h_{k+1} = h_k + 2g_k$ and
$g_{k+1} = h_{k+1} + g_k$ with initial conditions $h_1 = g_1 = 1$.
Standard techniques for solving linear recurrences yield

$$g_{k} = \frac{1}{2\sqrt{3}}\big((2+\sqrt{3})^k -
(2-\sqrt{3})^k\big) > (2+\sqrt{3})^{k-1}.$$

The claim now follows from  
$|\Omega_L| > |P_k| = g_n > (2+\sqrt{3})^{n-1}$.
\qedclaim

The three claims together with Lemma~\ref{lem:hour-glass} yield
$\Phi_{M_3(T_n)} \leq 1/(2+\sqrt{3})^{n-1}$. Which implies the theorem
via Fact~T.
\qed
\subsection{Slow mixing for 3-orientations with sub-linear maximum degree}

\begin{wrapfigure}[10]{r}{0.33\textwidth}
\vskip-6mm
\centering
\includegraphics[width=0.30\textwidth]{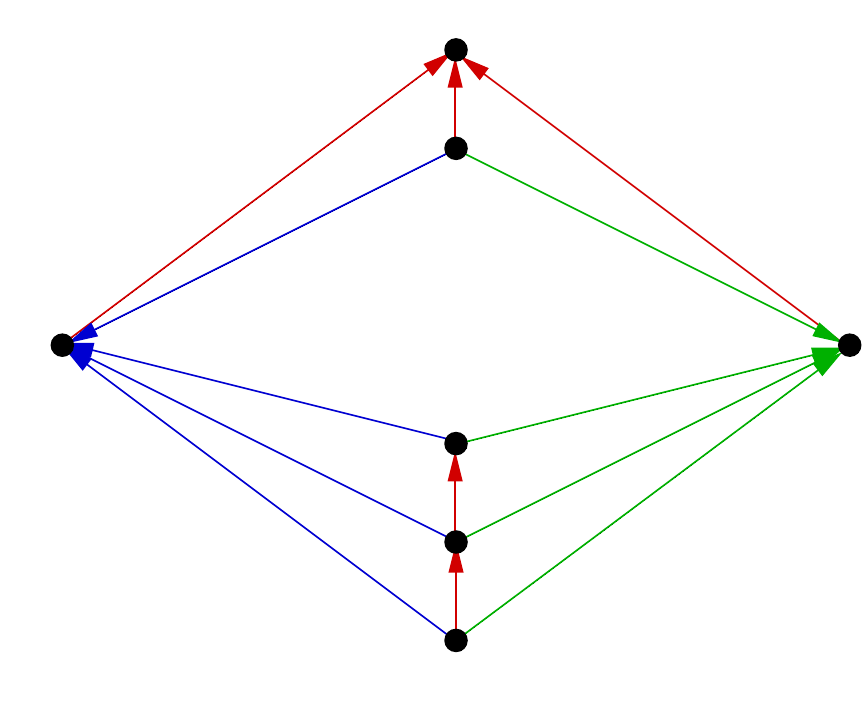}
\vskip-6mm
\caption{\label{fig:subdiv}
The gadget replacing $\{x_i,x_{i+1}\}$ in $T_n(m)$.}
\end{wrapfigure}
As announced we now modify $T_n$ to prove slow mixing for 
Schnyder woods of triangulations with a sub-linear maximum
degree. For a given $m\in\NN$ the triangulation $T_n(m)$
is constructed by replacing each edge $\{x_i,x_{i+1}\}$ with $i\geq 1$ 
by a path $x_i,y_{i,1},\dots,y_{i,m},x_{i+1}$. Each vertex 
$y_{i,j}$ is also made adjacent to to $v_i$ and $w_i$,
see Figure~\ref{fig:subdiv}. The resulting triangulation $T_n(m)$
has $3n+4 + (n-1)m$ vertices and its maximum degree is $\max\{
2n+3, m+5\}$. 


The definition of $\Omega_L,\Omega_c,\Omega_R$ for $T_n(m)$ is 
the same as for~$T_n$. This is again a hour glass partition, i.e.,
there is no direct transition between  $\Omega_L$ and $\Omega_R$.
Replacing the red edges $(x_i,x_{i+1})$ in Figure~\ref{fig:Tn}
by the colored gadget of Figure~\ref{fig:subdiv} yields the 
unique Schnyder wood $S$ of $\Omega_c$. To estimate $|\Omega_L|$ 
we again look at the set~$P_n$ of directed paths from $x_0$ to $x_n$
in $S$. Since there are $2m+3$ directed path from $x_i$ to $x_{i+1}$
we get $|P_n| > (2m+3)^{n-1}$. From $|\Omega_L| > (2m)^{n-1}$
we obtain:

\begin{proposition}\label{prop:degree}
  Let $T_n(2n)$ be the above triangulation on $2n^2 +n +4$ vertices with
  maximum degree $2n+5$.  The triangle flip Markov chain $M_3$ on
  3-orientations of the triangulation $T_n(n)$ has $\Tmix >
  2^{(n-1)\log(4n)-2}$.
\end{proposition}

\subsection{Slow mixing for $\alpha$-orientations with constant degree}

In~\cite{MRST-2016} and in this paper there are proofs for 
rapid mixing of the face flip Markov chain for $\alpha$-orientations
on graphs with small constant maximum degree and negative results in the sense
of slow mixing of these Markov chains for graphs with large maximum degree.
Could it be that the face flip Markov chain for $\alpha$-orientations
is rapidly mixing for all graphs of small maximum degree?
In this subsection we show that this is not the case.

\begin{figure}[htb]
    \centerline{\input{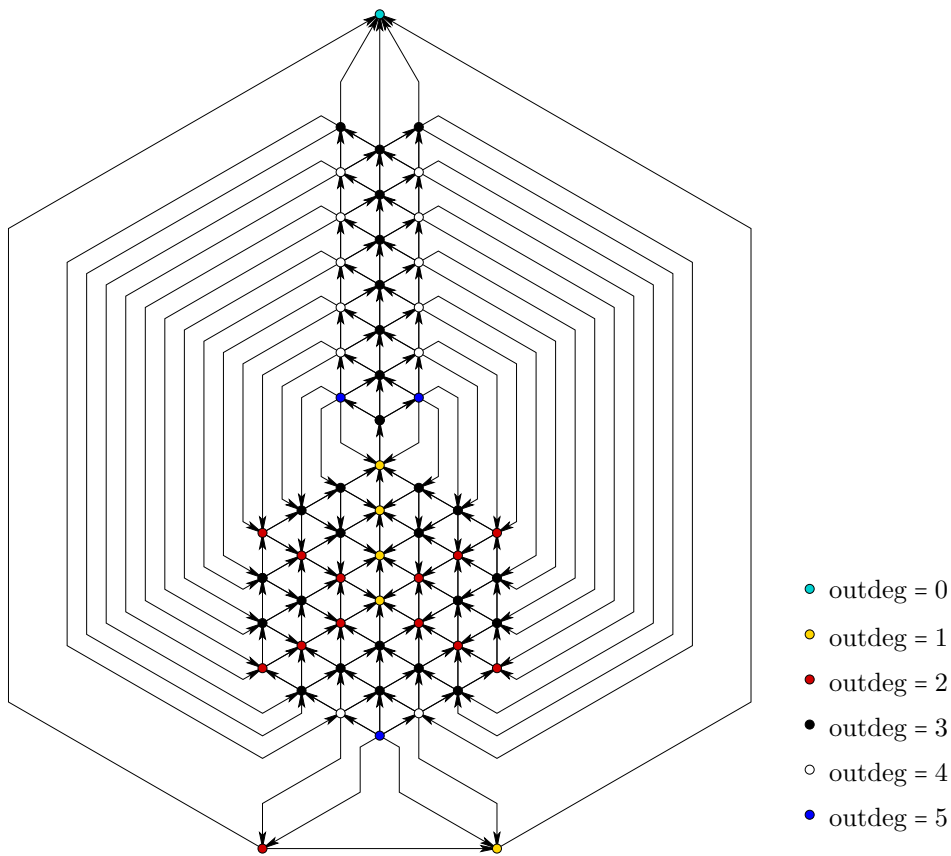}}
    \caption{The graph $G_3$, the degrees prescribed 
by $\alpha$, and the central $\alpha$-orientation.\label{fig:Gn}}
    \end{figure}

Our example family $G_k$ is obtained from $T_{3k-2}$. In $T_{3k-2}$ remove all
edges incident to $x_0$ except those connecting to $a_g$ and $a_b$. Let $H_k$
be a  patch taken from the triangular grid whose boundary is a regular hexagon
with side length $k$, i.e., each side has $k+1$ vertices, and in total $H_k$
has $3(k^2+k) + 1$ vertices. Now identify two opposite corners of $H_k$ with
the vertices $x_0$ and $x_1$ of~$T_{3k-2}$.  Label the vertices on the left
boundary of $H_k$ as $v'_0=x_1,v'_1,\ldots, v'_{3k-4},v'_{3k-3}=x_0$ and on
the right boundary as $w'_0=x_1,w'_1,\ldots, w'_{3k-4},w'_{3k-3}=x_0$. Add the
missing edges to make $v_i,v'_i,v'_{i+1}$ and $w_i,w'_i,w'_{i+1}$ triangles
for $i=1,\ldots 3k-2$. Finally add the edges from $v'_{3k-1}$ to~$a_b$ and
from $w'_{3k-1}$ to~$a_g$. Figure~\ref{fig:Gn} shows the result of the
construction for $k=3$.

   \calc_figscale{30}
    

The graph $G_k$ has $3(k^2 + 4k -1)$ vertices, the degrees are between 4 and
6. Let $\alpha$ be the function shown on the right part of
Figure~\ref{fig:Gn}, the values taken by $\alpha$ range from 0 to 5.
A key property of $\alpha$ is that except from the rigid edges which
connect to $a_g$ and $a_b$ there is exactly one out-edge of $H_k$, i.e.,
one edge directed from a vertex of $H_k$ to a vertex outside of $H_k$.

Let $\Omega$ denote the set of $\alpha$-orientation of $G_k$ We define
a hour glass partition $\Omega_L,\Omega_c,\Omega_R$ of this set. Let
$(y,z)$ be the unique non-rigid out-edge of $H_k$.  The edge $(y,z)$
is called {\em left} if $z\in\{v_1,\ldots,v_n\}$, it is {\em right} if
$z\in\{w_1,\ldots,w_n\}$ and it is {\em central} if $z=x_1$.  Now
$\Omega_L,\Omega_c,\Omega_R$ are the $\alpha$-orientations where the
edge $(y,z)$ is left, central, and right respectively.

It is clear that there is no transition between elements
of $\Omega_L$ and $\Omega_R$, i.e., the partition has the 
hour glass property. The set $\Omega_c$ has precisely one 
element, this is the orientation shown in
Figure~\ref{fig:Gn}. The size of $\Omega_L$ is at least as
large as the size of the set for $T_{3k-2}$ which has been shown to be
$> 3.73^{3k - 4}$. With the implied bound on the conductance and Fact~T
we obtain our last theorem.

\begin{theorem}\label{thm:lower-small-deg}
  The triangulation $G_k$ on $3(k^2 + 4k -1)$ vertices with maximum degree $6$
  and the function~$\alpha$ shown in Figure~\ref{fig:Gn} have 
  a slow mixing face flip Markov chain $M_\alpha$, more precisely 
  $\Tmix(M_\alpha) > 3.73^{3(k-2)}$.
\end{theorem}

\small
\advance\bibitemsep-3pt
\bibliography{newlit}

\newcommand{\etalchar}[1]{$^{#1}$}
\begin{thebibliography}{dFOdM01}

\bibitem[ABF{\etalchar{+}}13]{abfkku-ccoc-13}
Muhammad~Jawaherul Alam, Therese Biedl, Stefan Felsner, Michael Kaufmann,
  Stephen~G. Kobourov, and Torsten Ueckerdt.
\newblock Computing cartograms with optimal complexity.
\newblock {\em Discr. and Comput. Geom.}, 50:784--810, 2013.

\bibitem[BH12]{BH-2012}
Lali Barri{\`{e}}re and Clemens Huemer.
\newblock 4-labelings and grid embeddings of plane quadrangulations.
\newblock {\em Discr. Math.}, 312:1722--1731, 2012.

\bibitem[BSVV08]{BSVV-08}
Ivona Bez{\'{a}}kov{\'{a}}, Daniel Stefankovic, Vijay~V. Vazirani, and Eric
  Vigoda.
\newblock Accelerating simulated annealing for the permanent and combinatorial
  counting problems.
\newblock {\em SIAM J. Comput.}, 37:1429--1454, 2008.

\bibitem[Cre09]{Cr-09}
P{\'{a}}id{\'{\i}}~J. Creed.
\newblock Sampling {E}ulerian orientations of triangular lattice graphs.
\newblock {\em J. of Discr. Alg.}, 7:168--180, 2009.

\bibitem[dFOdM01]{dFOdM-tao-01}
Hubert de~Fraysseix and Patrice Ossona~de Mendez.
\newblock On topological aspects of orientations.
\newblock {\em Discr. Math.}, 229:57--72, 2001.

\bibitem[DG98]{DG-1998}
Martin Dyer and Catharine Greenhill.
\newblock A more rapidly mixing {M}arkov chain for graph colourings.
\newblock {\em Rand. Struct. and Alg.}, 13:285--317, 1998.

\bibitem[DSC93]{DSC-1993}
Persi Diaconis and Laurent Saloff-Coste.
\newblock Comparison theorems for reversible {M}arkov chains.
\newblock {\em An. of Appl. Prob.}, 3:696--730, 1993.

\bibitem[Fel04a]{f-gga-04}
S.~Felsner.
\newblock {\em Geometric Graphs and Arrangements}.
\newblock Vieweg Verlag, 2004.

\bibitem[Fel04b]{Fe-04}
Stefan Felsner.
\newblock Lattice structures from planar graphs.
\newblock {\em Electr. J. Combin.}, 11(1):24p., 2004.

\bibitem[FFNO11]{ffno-bsspbc-10}
Stefan Felsner, {\'E}ric Fusy, Marc Noy, and David Orden.
\newblock Bijections for {B}axter families and related objects.
\newblock {\em J. Combin. Theory Ser. A}, 18:993--1020, 2011.

\bibitem[FHKO10]{fhko-blpqr-10}
Stefan Felsner, Clemens Huemer, Sarah Kappes, and David Orden.
\newblock Binary labelings for plane quadrangulations and their relatives.
\newblock {\em Discr. Math. and Theor. Comp. Sci.}, 12:3:115--138, 2010.

\bibitem[FK09]{fk-uld-09}
Stefan Felsner and Kolja Knauer.
\newblock {ULD}-lattices and {$\Delta$}-bonds.
\newblock {\em Comb., Probab. and Comput.}, 18(5):707--724, 2009.

\bibitem[FR04]{FR-2004}
Johannes Fehrenbach and Ludger R\"uschendorf.
\newblock {M}arkov chain algorithms for {E}ulerian orientations and
  3-colourings of 2-dimensional cartesian grids.
\newblock {\em Statistics \& Decisions}, 22:109--130, 2004.

\bibitem[FZ08]{FZ-08}
Stefan Felsner and Florian Zickfeld.
\newblock On the number of planar orientations with prescribed degrees.
\newblock {\em Electr. J. Combin.}, 15:41p., 2008.

\bibitem[JSV04]{JSV-Permanent-04}
Mark Jerrum, Alistair Sinclair, and Eric Vigoda.
\newblock A polynomial-time approximation algorithm for the permanent of a
  matrix with nonnegative entries.
\newblock {\em J. ACM}, 51:671--697, 2004.

\bibitem[LPW09]{LevinPeresWilmer}
David Levin, Yuval Peres, and Elizabeth Wilmer.
\newblock {\em {M}arkov Chains and Mixing Times}.
\newblock AMS, 2009.

\bibitem[LRS95]{LuRaSi95}
M.~Luby, D.~Randall, and A.~Sinclair.
\newblock Markov chain algorithms for planar lattice structures.
\newblock In {\em Proc. FOCS}, pages 150--159, 1995.

\bibitem[MRST12]{MRST-2011}
Sarah Miracle, Dana Randall, Amanda~Pascoe Streib, and Prasad Tetali.
\newblock Mixing times of {M}arkov chains on 3-orientations of planar
  triangulations.
\newblock In {\em Proc. {A}of{A}'12}, Proc. AQ, pages 413--424. Discr. Math.
  and Theor. Comp. Sci., 2012.
\newblock full version arXiv:1202.4945.

\bibitem[MRST16]{MRST-2016}
Sarah Miracle, Dana Randall, Amanda~Pascoe Streib, and Prasad Tetali.
\newblock Sampling and counting 3-orientations of planar triangulations.
\newblock {\em SIAM J. Discr. Math.}, 30:801--831, 2016.

\bibitem[MW96]{MW-96}
Milena Mihail and Peter Winkler.
\newblock On the number of {E}ulerian orientations of a graph.
\newblock {\em Algorithmica}, 16:402--414, 1996.

\bibitem[Pro97]{Propp-DL-97}
James Propp.
\newblock Generating random elements of finite distributive lattices.
\newblock {\em Electr. J. Combin.}, 4:12 pp., 1997.

\bibitem[PS06]{ps-ocst-06}
D.~Poulalhon and G.~Schaeffer.
\newblock Optimal coding and sampling of triangulations.
\newblock {\em Algorithmica}, 46:505--527, 2006.

\bibitem[PW96]{PW-exact-96}
James~G. Propp and David~B. Wilson.
\newblock Exact sampling with coupled {M}arkov chains and applications to
  statistical mechanics.
\newblock {\em Rand. Struct. and Alg.}, 9(1\&2):223--252, 1996.

\bibitem[RT97]{RT-1997}
Dana Randall and Prasad Tetali.
\newblock Analyzing {G}lauber dynamics by comparison of {M}arkov chains.
\newblock {\em J. of Math. Phys.}, 41:1598--1615, 1997.

\bibitem[Sch89]{s-pgpd-89}
W.~Schnyder.
\newblock Planar graphs and poset dimension.
\newblock {\em Order}, 5:323--343, 1989.

\bibitem[SJ89]{SJ-89}
Alistair Sinclair and Mark Jerrum.
\newblock Approximate counting, uniform generation and rapidly mixing {M}arkov
  chains.
\newblock {\em Inf. Comput.}, 82:93--133, 1989.

\end{thebibliography}
\bibliographystyle{alpha}

\end{document}